\documentclass[letterpaper, 10 pt, conference]{ieeeconf}
\IEEEoverridecommandlockouts

\usepackage{amsmath}
\usepackage{graphicx}
\usepackage{cite}
\usepackage{authblk}
\usepackage{amssymb}

\begin{document}
\title{Decentralized Event-Triggered Consensus of Linear Multi-agent Systems under Directed Graphs}
\author{Eloy Garcia$^\ast$, Yongcan Cao, Xiaofeng Wang, and David W. Casbeer
\thanks{This work has been supported in part by AFOSR LRIR No. 12RB07COR.}
\thanks{A preliminary version of this manuscript has been submitted to the 2015 American Control Conference.}
\thanks{E. Garcia is a contractor (Infoscitex Corp.) with the Control Science Center of Excellence, Air Force Research Laboratory, Wright-Patterson AFB, OH 45433.}
\thanks{Y. Cao and D. Casbeer are with the Control Science Center of Excellence, Air Force Research Laboratory, Wright-Patterson AFB, OH 45433.}
\thanks{X. Wang is with the Department of Electrical Engineering, University of South Carolina, Columbia, SC, 29208.}
\thanks{$^\ast$ Corresponding Author: \ttfamily{elgarcia@infoscitex.com}}
}

\newtheorem{theorem}{Theorem}
\newtheorem{lemma}{Lemma}
\newtheorem{remark}{Remark}

\maketitle 
\begin{abstract}
An event-triggered control technique for consensus of multi-agent systems with general linear dynamics is presented. This paper extends previous work to consider agents that are connected using directed graphs. Additionally, the approach shown here provides asymptotic consensus with guaranteed positive inter-event time intervals. This event-triggered control method is also used in the case where communication delays are present. For the communication delay case we also show that the agents achieve consensus asymptotically and that, for every agent, the time intervals between consecutive transmissions is lower-bounded by a positive constant.

\end{abstract}

\section{Introduction} \label{sec:intro}
Cooperative control of multi-agent systems is an active research area with broad and relevant applications in commercial, academic and military areas~\cite{RenBeard07}. The design of decentralized and scalable control algorithms provides the necessary coordination for a group of agents to outperform a single or a number of systems operating independently. In general, agents use a limited bandwidth communication channel to broadcast information. Thus, continuous communication among agents is not possible to implement. Further, periodic communication schemes require global synchronization of sample periods and broadcasting time instants which are difficult to achieve in a decentralized setting. On the other hand, event-based communication offers a highly decentralized way to determine broadcasting time instants, that is, each agent is able to decide when to transmit measurements based only on locally available information.

\addtolength{\abovedisplayskip}{-1mm}
\addtolength{\belowdisplayskip}{-1mm}
\setlength{\textfloatsep}{2.0pt plus 1.0pt minus 1.0pt}
\setlength{\floatsep}{2.0pt plus 1.0pt minus 1.0pt}

In the present paper we address the event-triggered consensus problem where agents are described by general linear dynamics and are connected using directed graphs. In addition, we consider the case where communication among agents is subject to communication delays. Different from periodic (or time-triggered) implementations, in the context of event-triggered control, information or measurements are not transmitted periodically in time but they are triggered by the occurrence of certain events. In event-triggered broadcasting  \cite{Astrom08}, \cite{Donkers10}, \cite{DePersis13}, \cite{Garcia13}, and \cite{Tabuada07}, a subsystem sends its local state to the network only when it is necessary, that is, only when a measure of the local subsystem state error is above a specified threshold. Event-triggered control strategies have been used for stabilization of multiple coupled subsystems as in \cite{Garcia12}, \cite{Stocker13}, and \cite{Guinaldo11}. Consensus problems have also been studied using these techniques \cite{Dimarogonas12}, \cite{Garcia13b}, \cite{Seyboth13}, \cite{Yu12}, \cite{YinYue13b}, \cite{ChenHao12}, \cite{GuoDimarogonas13}. Event-triggered control provides a more robust and efficient use of network bandwidth. Its implementation in multi-agent systems also provides a highly decentralized way to schedule transmission instants which does not require synchronization compared to periodic sampled-data approaches. 

Consensus problems where all agents are described by general linear models have been considered by different authors \cite{LiDuan10}, \cite{LiDuan11}, \cite{MaZhang10}, \cite{Ren08}, \cite{Scardovi09}, \cite{SuHuang12}, \cite{Tuna08}, and \cite{Tuna09}. In these papers it is assumed that continuous communication between agents is possible. The work in \cite{WenDuan13} considers the consensus problem of agents with linear dynamics under communication constraints. Specifically, the authors consider the existence of continuous communication among agents for finite intervals of time and the total absence of communication among agents for other time intervals, and the minimum rate of continuous communication to no communication is given.

Event-triggered consensus of agents with linear dynamics and limited communication was recently explored in \cite{LiuHill12} and \cite{Zhu14}. In our previous work \cite{Garcia14auto} we proposed a novel approach in which each agent implements models of the decoupled dynamics of each one of its neighbors and uses the model states to compute the local control input. This approach offered better performance than the Zero-Order-Hold (ZOH) approach used in \cite{LiuHill12} and \cite{Zhu14} where the updates from neighbors are kept constant by the local agent. A similar model-based framework was proposed in \cite{Demir12b} where only constant thresholds were used. 
One of the main limitations of the ZOH approach \cite{LiuHill12}, \cite{Zhu14} is that it is not capable to keep up with unstable trajectories and updates need to be generated more frequently. In consensus with general linear dynamics, unstable systems are one of the most interesting cases to analyze. On the other hand the model-based approach in \cite{Garcia14auto} provides better estimates of neighbors and reduces generation of events as agents converge to similar unstable trajectories.  

The present paper describes a method for designing event thresholds that offers two important advantages with respect to \cite{Garcia14auto}. First, the results in this paper extend the work in \cite{Garcia14auto} to consider linear multi-agent systems that are inter-connected using directed graphs in contrast to the less general case studied in those papers where only undirected graphs were considered. Second, the event-triggered control strategy in this paper provides asymptotic consensus, while guaranteeing positive inter-event time intervals, compared to the results in \cite{Garcia14auto} where the difference between any two states can only be bounded but asymptotic convergence is not guaranteed. In addition, we extend the event-triggered consensus approach proposed here to consider the case where transmission of information among agents is subject to communication delays. Concerning the event-based consensus problem of linear systems, the paper \cite{Garcia14CDC} considered communication delays but, similar to \cite{Garcia14auto}, the results only addressed undirected graphs and only bounded consensus could be obtained.

The remainder of this paper is organized as follows. Section \ref{sec:Problem} provides a short background on graph theory and describes the problem. An event-triggered control strategy that achieves asymptotic consensus of multi-agent systems which are represented by general linear dynamics and connected by means of directed graphs is presented in Section \ref{sec:cons}. Section \ref{sec:Delays} provides similar results for the case of communication delays. An illustrative example is shown in Section \ref{sec:Examples}. Section \ref{sec:Conclusions} concludes the paper.

\section{Preliminaries} \label{sec:Problem}

\textit{Notation}.
Let $I_n$ represent an identity matrix of size $n$. The notations $1_n$ and $0_n$ represent column vectors of all ones and all zeros, respectively. $\mathbb{R}$ and $\mathbb{C}$ denote the set of real numbers and the set of complex numbers, respectively. For any $s\in\mathbb{C}$, $Re(s)$ represents the real part of $s$. The symbol $\otimes$ denotes the Kronecker product. $J_\mu^{\lambda}$ represents a Jordan block of size $\mu$ corresponding to eigenvalue $\lambda$. 

\subsection{Graph Theory}
For a team of $n$ agents, the communication among them can be described by a directed graph $\mathcal{G}=\left\{\mathcal{V,E}\right\}$, where $\mathcal{V}=\left\{1,\ldots,N\right\}$ denotes the agent set and $\mathcal{E}\subseteq\mathcal{V}\times\mathcal{V}$ denotes the edge set. An edge $(i,j)$ in the set $\mathcal{E}$ denotes that agent $j$ can obtain information from agent $i$, but not necessarily \emph{vice versa}. For an edge $(i,j)\in \mathcal{E}$, agent $i$ is a neighbor of agent $j$. The set $\mathcal{N}_j$ is called the set of neighbors of agent $j$, and $N_j$ is its cardinality.
A directed path from agent $i$ to agent $j$ is a sequence of edges in a directed graph of the form $(i,p_1),(p_1,p_2),\ldots,(p_{\kappa-1},p_{\kappa})(p_\kappa,j)$, where $p_\ell\in\mathcal{V},~\forall \ell=1,\cdots,\kappa$. A directed graph is \emph{strongly connected} if there is a directed path from every agent to every other agent. A directed graph has \emph{a directed spanning tree} if there exists at least one agent with directed paths to all other agents.

The adjacency matrix $\mathcal{A}\in\mathbb{R}^{n\times n}$ of a directed graph $\mathcal{G}$ is defined by $a_{ij}=1$ if $(j,i)\in\mathcal{E}$ and $a_{ij}=0$ otherwise. The Laplacian matrix $\mathcal{L}$ of $\mathcal{G}$ is defined as $\mathcal{L=D-A}$, where $\mathcal{D}$ represents the degree matrix which is a diagonal matrix with entries $d_{ii}=\sum_{j\in\mathcal{N}_i}a_{ij}$. If a directed graph has a directed spanning tree, then the corresponding Laplacian matrix has only one eigenvalue equal to zero, $\lambda_1=0$, and the following holds for the remaining eigenvalues: $Re\left\{\lambda_i\right\}>0$, for $i=2,...,N$.

\subsection{Problem Statement}
We consider the consensus problem with agents described by linear dynamics and with limited communication constraints where information from neighbors is not available continuously but only at some time instants.
Event-triggered control implementations typically use a ZOH \cite{Tabuada07} to compute the control input and the state error in problems where continuous feedback is not available. Model-based approaches have been used more recently and it has been shown that they offer better performance by providing an estimate of the real state of a system between update intervals \cite{Demir12a}, \cite{Garcia13}.
The model-based approach generalizes the traditional ZOH event-triggered control strategy. In ZOH strategies the agents that receive information from agent $i$ maintain a piece-wise constant model of the state $x_i(t)$. The ZOH case is equivalent to implementing models when $A=0$ in \eqref{eq:models} below. However, the choice of ZOH is not suitable when considering general linear dynamics as it was in the case of single integrators \cite{Dimarogonas12}, \cite{Garcia13b}. Since trajectories can be unstable in general, a ZOH is not able to reduce communication as trajectories grow. In this case sensors need to generate events more frequently since the errors grow very quickly after each update. This situation increases communication and Zeno behavior may not be avoided. 
In contrast, the models are able to produce better estimates of real states than the ZOH and it is possible to show that Zeno behavior does not occur. Note that, in contrast to \cite{Dimarogonas12}, the focus of this work is in reducing the number of transmissions instead of reducing actuation updates as it was discussed in that reference. 

Consider a group of $N$ agents with fixed and directed communication graphs and fixed weights. Each agent can be described by the following:
\begin{align}
	\dot{x}_i(t)=Ax_i(t)+Bu_i(t),\ \ i=1,...,N \label{eq:agents}
\end{align}
 with
\begin{align}
	u_i(t)=cF\sum_{j\in\mathcal{N}_i}(y_i(t)-y_j(t)),\ \ i=1,...,N \label{eq:inputs}
\end{align}
where $x_i\in \mathbb{R}^n$, $u_i\in \mathbb{R}^m$. The variables $y_i\in \mathbb{R}^n$ represent a model of the $i^{th}$ agent's state using the decoupled dynamics:
\begin{align}
  	\dot{y}_i(t)&=Ay_i(t), \;t\in[t_{k_i},t_{k_i+1})  \label{eq:models}  \\
		{y}_i(t_{k_i})&=x_i(t_{k_i}),  \nonumber
\end{align}
for $i=1,...,N$. Define the local errors $e_i(t)=y_i(t)-x_i(t)$.
Every agent in the network implements a model of itself $y_i(t)$ and also models of its neighbors $y_j(t)$. Local events for agent $i$ are triggered by the occurrence of the event 
\begin{align}
	 \left\|e_i(t)\right\|= \beta \textbf{e}^{-\lambda t}.   \label{eq:thre}
\end{align}
When agent $i$ triggers an event at time $t_{k_i}$, it will transmit its current state $x_i(t_{k_i})$ to its neighbors and agent $i$ and its neighbors will update their local models $y_i(t)$. Since agent $i$ and its neighbors use the same measurements to update the models and the model dynamics~\eqref{eq:models} represent the decoupled dynamics where all agents use the same state matrix, then the model states $y_i(t)$ implemented by agent $i$ and by its neighbors are the same. The model update process is similar for all agents $i=1,...,N$. In the presence of communication delays the previous statement will not hold and we will differentiate between $y_{i}(t)$, the model state of agent $i$ as seen by agent $i$, or the model with no delays; and $y_{i}^d(t)$, the model state of agent $i$ as seen by agents $j$, such that $i\in \mathcal{N}_j$, or the delayed model. More details concerning communication delays are presented in Section \ref{sec:Delays}. 

The local control input~\eqref{eq:inputs} is decentralized since it only depends on local information, that is, on the model states of the local agent and its neighbors.
Note that the difference between the agent dynamics~\eqref{eq:agents} and our proposed models~\eqref{eq:models} is given by the input term in~\eqref{eq:agents} and this input decreases as the agents approach a consensus state. It can also be seen that in the particular case when systems~\eqref{eq:agents} represent single integrator dynamics, then our models degenerate to ZOH models as in \cite{Dimarogonas12}, \cite{Garcia13b}.

\section{Event-triggered consensus with directed graphs} \label{sec:cons}
Let us start by defining the vectors $
x(t)=\begin{bmatrix}
	x_1(t)^T\ldots x_N(t)^T
\end{bmatrix}^T$,  $
y(t)=\begin{bmatrix}
	y_1(t)^T\ldots y_N(t)^T
\end{bmatrix}^T$, and $
e(t)=\begin{bmatrix}
	e_1(t)^T\ldots e_N(t)^T
\end{bmatrix}^T$. Then, the dynamics of the overall system can be written as follows:
\begin{align}
	  \dot{x}=\bar{A}x+\bar{B}y=(\bar{A}+\bar{B})x+\bar{B}e.  \label{eq:System}
\end{align}
where $\bar{A}=I_N\otimes A$, $\bar{B}=c\mathcal{L}\otimes BF$. Assume that the pair ($A,B$) is controllable. Then, for $\alpha>0$ there exists a (independent of the communication graph) symmetric and positive definite solution $P$ to
\begin{align}
	PA+A^TP-2PBB^TP+2\alpha P<0. \label{eq:LMI}
\end{align}
Let  
\begin{align}
  F&=-B^TP \label{eq:F}\\
	c&\geq1/Re(\lambda_2). \label{eq:c}
\end{align}
By selection of these controller gains we have that the matrix $\hat{A}$, defined in the following theorem, is a Hurwitz matrix.

Also, there exists a similarity transformation $S_L$ such that $\mathcal{L}_J=S^{-1}_L\mathcal{L}S_L$ is in Jordan canonical form. Define $S=S_L\otimes I_n$ and $\hat{x}=S^{-1}x$. Thus, we can obtain the transformed system dynamics
\begin{align}
\left.
	\begin{array}{l l}
 \dot{\hat{x}}\!\!\!&=S^{-1}\dot{x}  \\
    	     &=S^{-1}(\bar{A}+\bar{B})x+S^{-1}\bar{B}e  \\
					&=(\bar{A}+c\mathcal{L}_J\otimes BF)\hat{x}+(c\mathcal{L}_J\otimes BF)S^{-1}e
			\end{array}  \label{eq:Tranfdyn}  \right.
\end{align}
Since $\lambda_1(\mathcal{L})=0$ we have 
\begin{align}
	\mathcal{L}_J=\begin{bmatrix}
    0 & 0^T_{N-1} \\
	 0_{N-1} & J_{2:N}
\end{bmatrix}.  \nonumber
\end{align}
where the matrix $J_{2:N}\in \mathbb{C}^{(N-1)\times (N-1)}$ contains Jordan blocks corresponding to the eigenvalues $\lambda_2(\mathcal{L}),...,\lambda_N(\mathcal{L})$.

\begin{theorem}\label{th:consensus}
 Assume that the communication graph has a spanning tree and that the pair ($A,B$) is controllable. Define $F$ and $c$ as in~\eqref{eq:F} and~\eqref{eq:c}. Then agents~\eqref{eq:agents} with decentralized control inputs~\eqref{eq:inputs} based on models \eqref{eq:models} achieve consensus asymptotically when the local thresholds are defined as in \eqref{eq:thre}, where $\beta>0$ and $0<\lambda<\hat{\lambda}$. The parameter $\hat{\lambda}$ is such that $ \left\|\textbf{e}^{\hat{A}t}\right\|\leq\hat{\beta}\textbf{e}^{-\hat{\lambda} t}$, for $\hat{\beta}>0$ where $\hat{A}=\bar{A}_{2:N}+cJ_{2:N}\otimes BF$.
Furthermore, the agents do not exhibit Zeno behavior and the inter-event times $t_{k_i+1}-t_{k_i}$ for every agent $i=1,...,N$ are bounded by the \textit{positive} time $\tau$, that is
\begin{align}
  \tau\leq t_{k_i+1}-t_{k_i}  \label{eq:tk}
\end{align}	
where 
\begin{align}
\tau=\frac{\ln(1+\beta/K_3)}{\left\|A\right\|+\hat{\lambda}}  \label{eq:tausol}
\end{align}
and the positive parameter $K_3$ is defined in \eqref{eq:K3} below.
\end{theorem}

\textit{Proof}. 
Note that because of threshold \eqref{eq:thre}, the error $e_i$ is reset to zero at the event instants $t_{k_i}$, that is, $e_i(t_{k_i})=0$. Thus, the error $e_i$ satisfies $\left\|e_i(t)\right\|\leq \beta \textbf{e}^{-\lambda t}$ and we have that $ \left\|e(t)\right\|\leq\sqrt{N} \beta \textbf{e}^{-\lambda t}$.

Note that
\begin{align}
	\left.
	\begin{array}{l l}
	(c\mathcal{L}_J\otimes BF)S^{-1}&=(c\mathcal{L}_J\otimes BF)(S_L^{-1}\otimes I_n)  \\
	  &=\mathcal{L}_JS_L^{-1}\otimes cBF \\
	 &=\begin{bmatrix}
    0_{N}^T  \\
	 \Delta
\end{bmatrix}\otimes cBF. 
			\end{array}   \right. \nonumber
\end{align}
where $\Delta\in \mathbb{C}^{(N-1)\times N}$ is given by $\Delta=[0_{N-1} \ J_{2:N}]S_L^{-1}$. Therefore, we have that the transformed dynamics can be written as follows
\begin{align}
	  \dot{\hat{x}}_1=A\hat{x}_1  \label{eq:x1hat}
\end{align}
and
\begin{align}
	  \dot{\hat{x}}_{2:N}=\hat{A}\hat{x}_{2:N}+\hat{B}e  \label{eq:x2Nhat}
\end{align}
where $\hat{B}=c\Delta\otimes BF$. Matrix $\hat{A}$ is Hurwitz, then, there exist $\hat{\beta}$ and $\hat{\lambda}$ both greater than zero such that $\left\|\textbf{e}^{\hat{A}t}\right\|\leq \hat{\beta}\textbf{e}^{-\hat{\lambda} t}$.

The consensus problem has been transformed into the stabilization problem of system \eqref{eq:x2Nhat}. The response of \eqref{eq:x2Nhat} can be bounded as follows
\begin{align}
	\left.
	\begin{array}{l l}
	 \left\|\hat{x}_{2:N}(t)\right\|  \\
	=\left\|\textbf{e}^{\hat{A}t}\hat{x}_{2:N}(0) + \int_0^t \textbf{e}^{\hat{A}(t-s)}\hat{B}e(s) ds \right\|  \\
 \leq \hat{\beta}\hat{x}_0\textbf{e}^{-\hat{\lambda} t} + \int_0^t \hat{\beta}\textbf{e}^{-\hat{\lambda}(t-s)}\left\|\hat{B}\right\| \sqrt{N}\beta \textbf{e}^{-\lambda s}ds \\
\leq \hat{\beta}\hat{x}_0\textbf{e}^{-\hat{\lambda} t} + \frac{\sqrt{N}\beta\hat{\beta}\left\|\hat{B}\right\|}{\hat{\lambda}-\lambda}\big(\textbf{e}^{-\lambda t}-\textbf{e}^{-\hat{\lambda} t}\big)
  \end{array}   \right. \label{eq:xhat}
\end{align}
where, by abuse of notation, we denote $\hat{x}_0=\left\|\hat{x}_{2:N}(0)\right\|$. Note that 
\begin{align}
	   \lim_{t\rightarrow \infty}\left\|\hat{x}_{2:N}(t)\right\| = 0  \label{eq:xhatinf}
\end{align}
that is, the transformed states $\hat{x}_{2:N}(t)$ are asymptotically stable.

In order to show that the same condition guarantees asymptotic consensus we use the similarity transformation $S$. Note that $\lim_{t\rightarrow \infty}\hat{x}(t)=\begin{bmatrix}\lim_{t\rightarrow \infty}\hat{x}_1(t)^T\ 0\; ...\;0\end{bmatrix}^T$. 

Use the transformation $S=S_L\otimes I_n$ to obtain the original state $x$ from the states $\hat{x}$. Note that the first column of $S_L$ contains the right eigenvector of $\mathcal{L}$ associated with $\lambda_1=0$. Let $S_{L^{2:N}}$ denote the remaining columns of $S_L$, then we can write
\begin{align}
	 \lim_{t\rightarrow \infty}x(t)&=S\lim_{t\rightarrow \infty}\hat{x}(t)
		=\alpha \begin{bmatrix}
		 \lim_{t\rightarrow \infty}\hat{x}_1(t)\\
		\lim_{t\rightarrow \infty}\hat{x}_1(t)\\
		\vdots \\
		\lim_{t\rightarrow \infty}\hat{x}_1(t)
		\end{bmatrix}  \label{eq:Transfchat}
\end{align}
and the agents achieve consensus asymptotically. 

Note that for given controller parameters \eqref{eq:F} and \eqref{eq:c} the convergence rate of the event-triggered consensus algorithm is proportional to the selection of parameters $\beta$ and $\lambda$ as it can be seen in \eqref{eq:xhat}. 

In order to establish a positive lower-bound on the inter-event times (as a function of the selected convergence rate parameters $\beta$ and $\lambda$) for each agent $i=1,...,N$, we study the dynamics of the errors $e_i$, $i=1,...,N$.  
\begin{align}
	   \dot{e}_i=\dot{y}_i-\dot{x}_i=Ae_i-Bu_i=Ae_i-cBFz_i  \label{eq:eidot}
\end{align}
for $t\in[t_{k_i},t_{k_i+1})$, where
\begin{align}
   z_i=\sum_{j\in\mathcal{N}_i}(y_i(t)-y_j(t)). 
\end{align}
For the term $z_i$ the following holds
\begin{align}
\left.
	\begin{array}{l l}
   \left\|z_i\right\|\!\!\!&\leq \left\|z\right\|  \\
	  &\leq \left\|\mathcal{L}_n x\right\| + \left\|\mathcal{L}_n\right\|\left\|e\right\|
			\end{array}  \label{eq:zibound}  \right.
\end{align}
where $\mathcal{L}_n=\mathcal{L}\otimes I_n$. Also, we have that
\begin{align}
\left.
	\begin{array}{l l}
  \mathcal{L}_n x\!\!\!&=\mathcal{L}_n S\hat{x}=(\mathcal{L}\otimes I_n)(S_L\otimes I_n)\hat{x} \\
	 &=(\mathcal{L}S_L\otimes I_n)\hat{x}=(S_L\mathcal{L}_J\otimes I_n)\hat{x}
			\end{array}  \label{eq:Lnx}  \right.
\end{align}
Note that
\begin{align}
	\left.
	\begin{array}{l l}
	S_L\mathcal{L}_J\!\!\!&=S_L\begin{bmatrix}
    0 & 0_{N-1}^T \\
	  0_{N-1} & J_{2:N}
\end{bmatrix}= \begin{bmatrix}
    0_N & \Theta
\end{bmatrix} 
			\end{array}   \right. \nonumber
\end{align}
where $\Theta\in \mathbb{C}^{N\times (N-1)}$ is given by $\Theta=S_L\begin{bmatrix}
    0_{N-1}^T \\ J_{2:N}
\end{bmatrix} $. Then, we can write the following
\begin{align}
\left.
	\begin{array}{l l}
  \left\|\mathcal{L}_n x\right\|\!\!\!&=\left\|(\begin{bmatrix}
    0_N & \Theta
\end{bmatrix}\otimes I_n) \hat{x}\right\| \\
&=\left\|(\Theta\otimes I_n)\hat{x}_{2:N}\right\| \\
   &\leq \hat{\Theta} \left\|\hat{x}_{2:N}\right\| 
			\end{array}  \label{eq:LnxN}  \right.
\end{align}
where $\hat{\Theta}=\left\|\Theta\right\|$. From \eqref{eq:eidot}, \eqref{eq:zibound}, and \eqref{eq:LnxN} we obtain
\begin{align}
\left.
	\begin{array}{l l}
   \frac{d}{dt}\left\|e_i\right\|\!\!\!&\leq \left\|A\right\|\left\|e_i\right\| +\left\|cBF\right\| (\left\|\mathcal{L}_n x\right\| + \left\|\mathcal{L}\right\|\left\|e\right\|)  \\
	  &\leq  \left\|A\right\|\left\|e_i\right\| + \left\|cBF\right\|\Big(\hat{\Theta}\big(\hat{\beta}\hat{x}_0\textbf{e}^{-\hat{\lambda} t}  \\ 
		&~~+ \frac{\sqrt{N}\beta\hat{\beta}\left\|\hat{B}\right\|}{\hat{\lambda}-\lambda}(\textbf{e}^{-\lambda t}-\textbf{e}^{-\hat{\lambda} t})\big) \\
		&~~+ \left\|\mathcal{L}\right\|\sqrt{N}\beta \textbf{e}^{-\lambda t} \Big)
			\end{array}  \label{eq:eidotb}  \right.
\end{align}
for $t\in[t_{k_i},t_{k_i+1})$, with $e_i(t_{k_i})=0$. The error response during the time interval $t\in[t_{k_i},t_{k_i+1})$ can be bounded as follows
\begin{align}
\left.
	\begin{array}{l l}
   \left\|e_i(t)\right\|\!\!\!&\leq \int_{t_{k_i}}^t \textbf{e}^{\left\|A\right\|(t-s)}\left\|cBF\right\|\big(K_1\textbf{e}^{-\hat{\lambda}s}+K_2\textbf{e}^{-\lambda s}\big)ds  \\
	&\leq \big(\textbf{e}^{\left\|A\right\|\tau}-\textbf{e}^{-\hat{\lambda}\tau}\big)\frac{K_1\left\|cBF\right\|}{\left\|A\right\|+\hat{\lambda}}\textbf{e}^{-\hat{\lambda}t_{k_i}} \\
	&~~+\big(\textbf{e}^{\left\|A\right\|\tau}-\textbf{e}^{-\lambda\tau}\big)\frac{K_2\left\|cBF\right\|}{\left\|A\right\|+\lambda}\textbf{e}^{-\lambda t_{k_i}}
			\end{array}  \label{eq:eib}  \right.
\end{align}
where $\tau=t-t_{k_i}$ and
\begin{align}
\left.
	\begin{array}{l l}
   K_1= \hat{\Theta}\Big(\hat{\beta}\hat{x}_0 - \frac{\sqrt{N}\beta\hat{\beta}\left\|\hat{B}\right\|}{\hat{\lambda}-\lambda}\Big)    \\ 
	K_2= \hat{\Theta}\frac{\sqrt{N}\beta\hat{\beta}\left\|\hat{B}\right\|}{\hat{\lambda}-\lambda}	+ \left\|\mathcal{L}\right\|\sqrt{N}\beta
			\end{array}  \nonumber  \right.
\end{align}
Thus, the time $\tau>0$ that it takes for the last expression in \eqref{eq:eib} to grow from zero, at time $t_{k_i}$, to reach the threshold $\beta\textbf{e}^{-\lambda t}=\beta\textbf{e}^{-\lambda(t_{k_i}+\tau)}$ is less or equal than the time it takes the error $\left\|e_i(t)\right\|$ to grow from zero, at time $t_{k_i}$, to reach the same threshold and generate the following event at time $t_{k_i+1}$, that is, $0<\tau\leq t_{k_i+1}-t_{k_i}$. Thus, we wish to find a lower-bound $\tau>0$ such that the following holds
\begin{align}
\left.
	\begin{array}{l l}
   &\big(\textbf{e}^{\left\|A\right\|\tau}-\textbf{e}^{-\hat{\lambda}\tau}\big)\frac{K_1\left\|cBF\right\|}{\left\|A\right\|+\hat{\lambda}}\textbf{e}^{-\hat{\lambda}t_{k_i}} \\
	&+\big(\textbf{e}^{\left\|A\right\|\tau}-\textbf{e}^{-\lambda\tau}\big)\frac{K_2\left\|cBF\right\|}{\left\|A\right\|+\lambda}\textbf{e}^{-\lambda t_{k_i}} \leq \beta\textbf{e}^{-\lambda(t_{k_i}+\tau)}
			\end{array}  \label{eq:taueq}  \right.
\end{align}
which can also be written as 
\begin{align}
\left.
	\begin{array}{l l}
   &\big(\textbf{e}^{\left\|A\right\|\tau}-\textbf{e}^{-\hat{\lambda}\tau}\big)\frac{K_1\left\|cBF\right\|}{\left\|A\right\|+\hat{\lambda}}\textbf{e}^{(\lambda-\hat{\lambda})t_{k_i}} \\
	&\ \ \ \ \ \ \ +\big(\textbf{e}^{\left\|A\right\|\tau}-\textbf{e}^{-\lambda\tau}\big)\frac{K_2\left\|cBF\right\|}{\left\|A\right\|+\lambda} \leq \beta\textbf{e}^{-\lambda \tau}
			\end{array}  \label{eq:tauU}  \right.
\end{align}
An explicit solution $\tau>0$ that guarantees \eqref{eq:tauU} can be found as follows. Let $\bar{\lambda}=\hat{\lambda}>\lambda$, 
then, the following two inequalities hold for any $\tau\geq 0$: 
\begin{align}
\left.
	\begin{array}{l l}
   &\big(\textbf{e}^{\left\|A\right\|\tau}-\textbf{e}^{-\hat{\lambda}\tau}\big)\frac{K_1\left\|cBF\right\|}{\left\|A\right\|+\hat{\lambda}}\textbf{e}^{(\lambda-\hat{\lambda})t_{k_i}} \\
	&~~+\big(\textbf{e}^{\left\|A\right\|\tau}-\textbf{e}^{-\lambda\tau}\big)\frac{K_2\left\|cBF\right\|}{\left\|A\right\|+\lambda} \\
	\leq \!\!\!&\big(\textbf{e}^{\left\|A\right\|\tau}-\textbf{e}^{-\hat{\lambda}\tau}\big)\frac{K_1\left\|cBF\right\|}{\left\|A\right\|+\hat{\lambda}}\textbf{e}^{(\lambda-\hat{\lambda})t_{k_i}} \\
	&~~+\big(\textbf{e}^{\left\|A\right\|\tau}-\textbf{e}^{-\bar{\lambda}\tau}\big)\frac{K_2\left\|cBF\right\|}{\left\|A\right\|+\lambda} 
 \end{array}  \label{eq:ineqtau1}  \right.
\end{align}
and
\begin{align}
\left.
	\begin{array}{l l}
\beta\textbf{e}^{-\bar{\lambda} \tau} \leq  \beta\textbf{e}^{-\lambda \tau}
			\end{array}  \label{eq:ineqtau2}  \right.
\end{align}
Then, the solution $\tau>0$ of 
\begin{align}
\left.
	\begin{array}{l l}
  &\big(\textbf{e}^{\left\|A\right\|\tau}-\textbf{e}^{-\hat{\lambda}\tau}\big)\frac{K_1\left\|cBF\right\|}{\left\|A\right\|+\hat{\lambda}}\textbf{e}^{(\lambda-\hat{\lambda})t_{k_i}} \\
	&~~+\big(\textbf{e}^{\left\|A\right\|\tau}-\textbf{e}^{-\bar{\lambda}\tau}\big)\frac{K_2\left\|cBF\right\|}{\left\|A\right\|+\lambda} = \beta\textbf{e}^{-\bar{\lambda} \tau}
 \end{array}  \label{eq:EqTau}  \right.
\end{align}
guarantees that inequality \eqref{eq:tauU} holds. Such solution is given by \eqref{eq:tausol}
where
\begin{align}
  K_3=\left\|cBF\right\|\Big(\frac{K_1\textbf{e}^{(\lambda-\hat{\lambda})t_{k_i}}}{\left\|A\right\|+\hat{\lambda}} +\frac{K_2}{\left\|A\right\|+\lambda} \Big).  \label{eq:K3}
\end{align}
By the selection $\hat{\lambda}>\lambda$, we have that $\textbf{e}^{(\lambda-\hat{\lambda})t_{k_i}}\leq 1$ for any $t_{k_i}\geq 0$, and the term $K_3$ remains bounded for any $t_{k_i}\geq 0$, ensuring that $\tau>0$. $\square$

\begin{remark} 
Note that the parameters $\beta$ and $\lambda$ do not need to be the same for all agents $i=1,...,N$. In general, each agent can use any $\beta_i>0$ and $0<\lambda_i<\hat{\lambda}$ and the consensus result follows by defining $\beta=\max_i\beta_i$ and $\lambda=\min_i\lambda_i$.
\end{remark}

\begin{remark}
It can be seen that if $\lambda>\hat{\lambda}$ then the second term in \eqref{eq:xhat} remains positive and asymptotic consensus is obtained. However, by making this selection, we try to impose a fast convergence of the state error with respect to the closed-loop response. By making this choice the inter-event time intervals will go to zero and continuous communication cannot be avoided. This can be clearly seen in the exponential term in \eqref{eq:K3} that will make $K_3$ to grow unbounded as time goes to infinity. 
\end{remark}

\begin{remark}
The parameters $\hat{\beta}$ and $\hat{\lambda}$ are related to the response of the closed-loop consensus protocol, which in turn is determined by the local matrices $A$, $B$, and $F$. These parameters also depend on the communication graph, in particular, on the second smallest eigenvalue of the Laplacian matrix $\lambda_2(\mathcal{L})$. As with many consensus algorithms, an estimate of the second smallest eigenvalue of the Laplacian matrix is required; this is the only global information needed by the agents.  Algorithms for distributed estimation of the second eigenvalue of the Laplacian have been presented in \cite{Aragues12}, \cite{Franceschelli09}. Readers are referred to these papers for details.
\end{remark}

\begin{remark} 
The convergence rate of the event-triggered algorithm is slower, as expected, compared to the case when continuous communication is possible. Such convergence rate is given only by the first term in \eqref{eq:xhat} where the parameters $\hat{\beta}$ and $\hat{\lambda}$ depend only on the system dynamics, the communication graph, and the chosen controller parameters. The selection of parameters $\beta$ and $\lambda$ provide a tradeoff between convergence rate and reduction of communication as measured by the minimum inter-event time intervals.
\end{remark}

\section{Event-triggered consensus with directed graphs and communication delays} \label{sec:Delays}
In this section we consider constant communication delays $d$.
Since the measurement updates will be delayed, the agents that receive information from agent $i$ will have a version of agent $i$'s model state that it is different than agent $i$'s version. Thus, it is necessary to distinguish between the model state as seen by the local agent and as seen by agents $j$, for $i\in\mathcal{N}_j$. Define the dynamics and update law of the model state of agent $i$ as seen by agent $i$ as
\begin{align}
  \dot{y}_{i}(t)=Ay_{i}(t),  \;  y_{i}(t_{k_i})=x_{i}(t_{k_i}). \label{eq:yii}
\end{align}
for $t\in[t_{k_i},t_{k_i+1})$. The measurement $x_{i}(t_{k_i})$ is transmitted by agent $i$ at time $t_{k_i}$ and will arrive to agents $j$, such that $i\in \mathcal{N}_j$, at time $t_{k_i}+d$. Let $y_i^d$ denote the state of the model of agent $i$ as seen by agents $j$ (the delayed model state of agent $i$), $i\in \mathcal{N}_j$. Define the dynamics and update law of $y_i^d$ as
\begin{align}
\begin{array}{l l}
  \dot{y}_i^d(t)=Ay_i^d(t),  \\
		y_i^d(t_{k_i}+d)=f_d(x_{i}(t_{k_i}),d) \label{eq:yij}
		\end{array}
\end{align}
for $t\in[t_{k_i}+d,t_{k_i+1}+d)$.
Since both, $y_{i}$ and $y_i^d$, use the same state matrix for their continuous evolution between their corresponding update instants, then we define
\begin{align}
   f_d(x_{i}(t_{k_i}),d)=e^{Ad}x_{i}(t_{k_i}). \label{eq:yij_prop}
\end{align}
In the presence of communication delays every agent $i=1,...,N$ will implement an additional model of itself. The first model is similar to the one used in the previous section and is represented by $y_i(t)$. The model state $y_i(t)$ is used by the local agent to compute the local error and to determine the local event time instants. The second model is represented by $y_i^d(t)$ which is a delayed model, equivalent to the models that other agents implement of agent $i$. The second model is updated at time instants $t_{k_i}+d$ using the update law in \eqref{eq:yij} and \eqref{eq:yij_prop}. The state of this model is used (along with model states $y_j^d$, $j\in \mathcal{N}_i$) to compute the local control inputs.

The agent dynamics are given by \eqref{eq:agents} and the control inputs are now defined as follows
\begin{align}
	u_i(t)=cF\sum_{j\in\mathcal{N}_i}(y_i^d(t)-y_j^d(t)),\ \ i=1,...,N. \label{eq:inputsd}
\end{align}
Define the state errors
\begin{align}
  e_{i}(t)&=y_{i}(t)-x_i(t), \label{eq:eii}  \\ 
  e_i^d(t)&=y_i^d(t)-x_i(t).  \label{eq:eij}
\end{align}
Note that $e_{i}(t_{k_i})=0$.
The dynamics of the overall system can be written as follows:
\begin{align}
	  \dot{x}=\bar{A}x+\bar{B}y_d=(\bar{A}+\bar{B})x+\bar{B}e_d  \label{eq:Systemd}
\end{align}
where $
y_d(t)=\begin{bmatrix}
	y_1^d(t)^T\ldots y_N^d(t)^T
\end{bmatrix}^T$, and $
e_d(t)=\begin{bmatrix}
	e_1^d(t)^T\ldots e_N^d(t)^T
\end{bmatrix}^T$. 

\begin{theorem}\label{th:consensusDel}
Assume that the communication graph has a spanning tree and that the pair ($A,B$) is controllable. Define $F$ and $c$ as in~\eqref{eq:F} and~\eqref{eq:c}. Then, there exists an $\epsilon>0$ such that for constant communication delays in the range $d\in[0,\epsilon)$ the agents~\eqref{eq:agents} with decentralized control inputs~\eqref{eq:inputsd} based on models \eqref{eq:models} achieve consensus asymptotically when the local thresholds are defined as in \eqref{eq:thre}
for $i=1,...,N$, where $\beta>0$ and $0<\lambda<\hat{\lambda}$. The parameter $\hat{\lambda}$ is such that $ \left\|\textbf{e}^{\hat{A}t}\right\|\leq\hat{\beta}\textbf{e}^{-\hat{\lambda} t}$, for $\hat{\beta}>0$ where $\hat{A}=\bar{A}_{2:N}+cJ_{2:N}\otimes BF$.
Furthermore, the agents do not exhibit Zeno behavior and the inter-event times $t_{k_i+1}-t_{k_i}$ for every agent $i=1,...,N$ are bounded by the \textit{positive} time $\tau$, that is
\begin{align}
  \tau\leq t_{k_i+1}-t_{k_i}  \nonumber
\end{align}	
where 
\begin{align}
\tau=\frac{\ln(1+\beta/H_3)}{\left\|A\right\|+\hat{\lambda}}   \label{eq:tausoldel}
\end{align}
and the positive parameter $H_3$ is defined in \eqref{eq:H3} below.
\end{theorem}
\textit{Proof}.
The error $e_{i}$ is used to trigger events. The time-dependent threshold is still given by \eqref{eq:thre}.
Then, there exist admissible delays $d\in[0,\epsilon)$ (where $\epsilon$ will be determined later in this proof) such that 
\begin{align}
\left\|e_i^d(t)\right\|\leq\gamma\textbf{e}^{-\lambda t}   \label{eq:thredel}
\end{align} 
for $i=1,...,N$, $\gamma>\beta$, and $t\in[t_{k_i},t_{k_i}+d)$. Note that because of the bound \eqref{eq:thredel} we have that $\left\|e_d\right\|\leq \sqrt{N}\gamma\textbf{e}^{-\lambda t}$.

Using the similarity transformation $S$ we can obtain the transformed system dynamics
\begin{align}
\left.
	\begin{array}{l l}
 \dot{\hat{x}}&=S^{-1}\dot{x}  \\
					&=(\bar{A}+c\mathcal{L}_J\otimes BF)\hat{x}+(c\mathcal{L}_J\otimes BF)S^{-1}e_d
			\end{array}  \label{eq:Tranfdynd}  \right.
\end{align}
Following similar steps as in Section \ref{sec:cons} we can write the transformed system dynamics as in \eqref{eq:x1hat} and
\begin{align}
	  \dot{\hat{x}}_{2:N}=\hat{A}\hat{x}_{2:N}+\hat{B}e_d.  \label{eq:x2Nhatd}
\end{align}
Recall that $\hat{A}$ is Hurwitz and $\left\|\textbf{e}^{\hat{A}t}\right\|\leq \hat{\beta}\textbf{e}^{-\hat{\lambda} t}$ holds for $\hat{\beta}>0$ and $\hat{\lambda}>0$.

The response of \eqref{eq:x2Nhatd} is now bounded as a function of the parameter $\gamma$ as follows
\begin{align}
	\left.
	\begin{array}{l l}
	 \left\|\hat{x}_{2:N}(t)\right\|  \leq \hat{\beta}\hat{x}_0\textbf{e}^{-\hat{\lambda} t} + \frac{\sqrt{N}\gamma\hat{\beta}\left\|\hat{B}\right\|}{\hat{\lambda}-\lambda}\big(\textbf{e}^{-\lambda t}-\textbf{e}^{-\hat{\lambda} t}\big)
  \end{array}   \right. \label{eq:xhatd}
\end{align}
where $0<\lambda<\hat{\lambda}$. Therefore the transformed states $\hat{x}_{2:N}(t)$ are asymptotically stable which means that asymptotic consensus is obtained for the real states $x(t)$ as it was shown in \eqref{eq:Transfchat}.

Note that once the local controller is specified (eqs. \eqref{eq:F} and \eqref{eq:c}) the convergence rate of the event-triggered consensus algorithm with communication delays is dictated by the selection of parameters $\beta$, $\gamma$, and $\lambda$.

Now, in order to determine the parameter $\epsilon>0$ (for given convergence rate parameters) such that for any constant delay in the range $d\in[0,\epsilon)$ the expression \eqref{eq:thredel} holds for a given delay $d\in[0,\epsilon)$, let us establish a bound on the response of the error $e_i^d(t)$. Here, we consider the general case where the delay can be larger than the inter-event time intervals. The error dynamics are given by
\begin{align}
	   \dot{e}_i^d=\dot{y}_i^d-\dot{x}_i=Ae_i^d-Bu_i=Ae_i^d-cBFz_i^d  \label{eq:eidotd}
\end{align}
for $t\in[t_{k_i}+d,t_{k_i+1}+d)$, where $z_i^d$ is given by 
\begin{align}
   z_i^d=\sum_{j\in\mathcal{N}_i}(y_i^d(t)-y_j^d(t)). 
\end{align}
For the term $z_i^d$ the following holds
\begin{align}
\left.
	\begin{array}{l l}
   \left\|z_i^d\right\|\!\!\!&\leq \left\|z\right\|  \\
	  &\leq \left\|\mathcal{L}_n x\right\| + \left\|\mathcal{L}_n\right\|\left\|e_d\right\|
			\end{array}  \label{eq:ziboundd}  \right.
\end{align}
Then, we can write the following
\begin{align}
\left.
	\begin{array}{l l}
  \left\|\mathcal{L}_n x\right\|\!\!\!  &\leq \hat{\Theta} \left\|\hat{x}_{2:N}\right\| 
			\end{array}  \nonumber  \right.
\end{align}
The error $e_i^d(t)$ is piece-wise continuous during the time interval $t\in[t_{k_i-1}+d,t_{k_i}+d)$. However, the response of the error $e_i^d(t)$ during that time interval depends on the value of the error $e_i(t^-_{k_i})$ at time $t^-_{k_i}$ (just before the event time at time $t_{k_i}$), as expected, since the delayed model state $y_i^d$ is updated using the update law \eqref{eq:yij}-\eqref{eq:yij_prop} with the information $x_i(t_{k_i})$ obtained at time $t_{k_i}$. 

\begin{figure}
	\begin{center}
		\includegraphics[width=8.4cm,height=6.5cm,trim=.8cm .1cm .8cm .4cm]{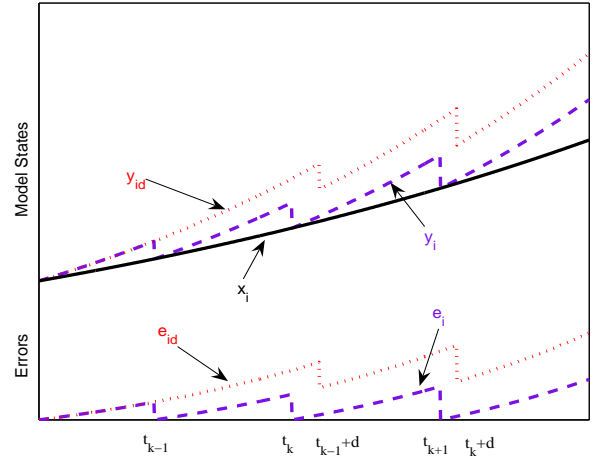}
	\caption{Relation between state $x_i$, model states $y_{i}$, $y_{i_d}$, and corresponding errors $e_{i}$, $e_{i_d}$.}
	\label{fig:FigureModErr}
	\end{center}
\end{figure}

Define the auxiliary model state variable $\tilde{y}_i^d(t)$ with dynamics given by
\begin{align}
\begin{array}{l l}
  \dot{\tilde{y}}_i^d(t)=A\tilde{y}_i^d(t),  \\
		\tilde{y}_i^d(t_{k_i})=x_{i}(t_{k_i}) \label{eq:auxyid}
		\end{array}
\end{align}
for $t\in[t_{k_i},t_{k_i}+d)$. Note that this model variable is not really implemented by the local agent but it is only used for the analysis to show convergence in the presence of delays. Because of the update law \eqref{eq:yij}-\eqref{eq:yij_prop} we have that $\tilde{y}_i^d(t)=y_i^d(t)$ for $t\in[t_{k_i-1}+d,t_{k_i}+d)$. Define the auxiliary error variable $\tilde{e}_i^d(t)=\tilde{y}_i^d(t)-x_i(t)$. Similarly, $\tilde{e}_i^d(t)=e_i^d(t)$ for $t\in[t_{k_i-1}+d,t_{k_i}+d)$. These relationships are illustrated in Fig. \ref{fig:FigureModErr}. Thus, by bounding the response of the error $\tilde{e}_i^d(t)$ for the time interval $t\in[t_{k_i},t_{k_i}+d)$ we are also establishing a bound on the response of the error ${e}_i^d(t)$  for the time interval $t\in[t_{k_i-1}+d,t_{k_i}+d)$.
 
The auxiliary error dynamics are bounded by
\begin{align}
\left.
	\begin{array}{l l}
   \frac{d}{dt}\left\|\tilde{e}_i^d\right\|\!\!\!&\leq  \left\|A\right\|\left\|\tilde{e}_i^d\right\| + \left\|cBF\right\|\Big(\hat{\Theta}\big(\hat{\beta}\hat{x}_0\textbf{e}^{-\hat{\lambda} t}  \\ 
		&~~+ \frac{\sqrt{N}\gamma\hat{\beta}\left\|\hat{B}\right\|}{\hat{\lambda}-\lambda}(\textbf{e}^{-\lambda t}-\textbf{e}^{-\hat{\lambda} t})\big) \\
		&~~+ \left\|\mathcal{L}\right\|\sqrt{N}\gamma \textbf{e}^{-\lambda t} \Big)
			\end{array}  \label{eq:eidotbd}  \right.
\end{align}
for $t\in[t_{k_i},t_{k_i}+d)$ with $\left\|\tilde{e}_i^d(t_{k_i})\right\|=\left\|e_i(t^-_{k_i})\right\|=\beta\textbf{e}^{-\lambda t_{k_i}}$. The initial condition in this expression represents the fact that a local event has been generated at time $t_{k_i}$.

\begin{figure}
	\begin{center}
		\includegraphics[width=8.4cm,height=6.5cm,trim=.8cm .2cm .8cm .4cm]{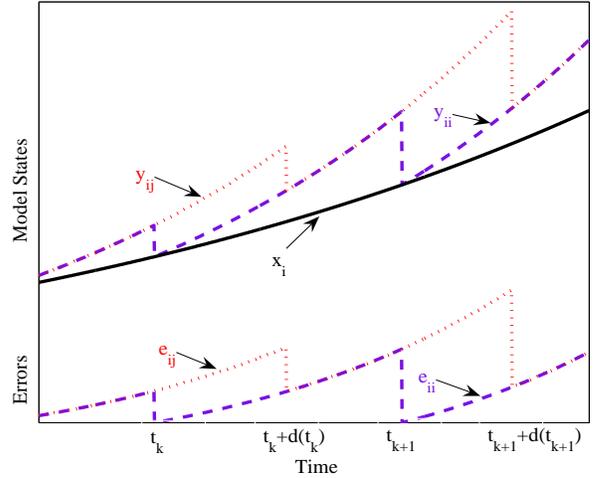}
	\caption{Relation between state $x_i$, model states $y_{i}$, $y_{i_d}$, and corresponding errors $e_{ii}$, $e_{i_d}$ when the delay is smaller than the inter-event time intervals.}
	\label{fig:FigureModErr2}
	\end{center}
\end{figure}

The auxiliary error response during the time interval $t\in[t_{k_i},t_{k_i}+d)$, with initial conditions $\left\|\tilde{e}_i^d(t_{k_i})\right\|=\beta\textbf{e}^{-\lambda t_{k_i}}$ can be bounded as follows
\begin{align}
\left.
	\begin{array}{l l}
   \left\|\tilde{e}_i^d(t)\right\|\!\!\!\!\!&\leq \beta\textbf{e}^{\left\|A\right\|(t-t_{k_i})}\textbf{e}^{-\lambda t_{k_i}}  \\
	  &~~+\int_{t_{k_i}}^t \textbf{e}^{\left\|A\right\|(t-s)}\left\|cBF\right\|\big(H_1\textbf{e}^{-\hat{\lambda}s}\!+\!H_2\textbf{e}^{-\lambda s}\big)ds  \\
	&\leq \beta\textbf{e}^{\left\|A\right\|d}\textbf{e}^{-\lambda t_{k_i}}  \\
	&~~+\big(\textbf{e}^{\left\|A\right\|d}-\textbf{e}^{-\hat{\lambda}d}\big)\frac{H_1\left\|cBF\right\|}{\left\|A\right\|+\hat{\lambda}}\textbf{e}^{-\hat{\lambda}t_{k_i}} \\
	&~~+\big(\textbf{e}^{\left\|A\right\|d}-\textbf{e}^{-\lambda d}\big)\frac{H_2\left\|cBF\right\|}{\left\|A\right\|+\lambda}\textbf{e}^{-\lambda t_{k_i}}
			\end{array}  \label{eq:auxeib}  \right.
\end{align}
where $d=t-t_{k_i}$ and
\begin{align}
\left.
	\begin{array}{l l}
   H_1= \hat{\Theta}\Big(\hat{\beta}\hat{x}_0 - \frac{\sqrt{N}\gamma\hat{\beta}\left\|\hat{B}\right\|}{\hat{\lambda}-\lambda}\Big)    \\ 
	H_2= \hat{\Theta}\frac{\sqrt{N}\gamma\hat{\beta}\left\|\hat{B}\right\|}{\hat{\lambda}-\lambda}	+ \left\|\mathcal{L}\right\|\sqrt{N}\gamma
			\end{array}  \nonumber  \right.
\end{align}
Note that if the delay is less than the minimum inter-event time intervals, there is no need for the auxiliary error variable since $e_i^d$ is continuous for $t\in[t_{k_i},t_{k_i}+d)$ and the same bound applies. This case is illustrated in Fig. \ref{fig:FigureModErr2}

In order to guarantee that \eqref{eq:thredel} holds for some desired convergence rate defined by the parameters $\lambda$ and $\gamma>\beta$ we need to find the range of values of $d$ such that the following inequality holds
\begin{align}
\left.
	\begin{array}{l l}
   \beta\textbf{e}^{\left\|A\right\|d}\textbf{e}^{-\lambda t_{k_i}} \\
	+\big(\textbf{e}^{\left\|A\right\|d}-\textbf{e}^{-\hat{\lambda}d}\big)\frac{H_1\left\|cBF\right\|}{\left\|A\right\|+\hat{\lambda}}\textbf{e}^{-\hat{\lambda}t_{k_i}} \\
	+\big(\textbf{e}^{\left\|A\right\|d}-\textbf{e}^{-\lambda d}\big)\frac{H_2\left\|cBF\right\|}{\left\|A\right\|+\lambda}\textbf{e}^{-\lambda t_{k_i}}  \leq \gamma\textbf{e}^{-\lambda(t_{k_i}+d)}
			\end{array}  \label{eq:dineq}  \right.
\end{align}
which can also be written as 
\begin{align}
\left.
	\begin{array}{r r}
   \beta\textbf{e}^{\left\|A\right\|d}\! + \!\big(\textbf{e}^{\left\|A\right\|d}\!-\!\textbf{e}^{-\hat{\lambda}d}\big)\frac{H_1\left\|cBF\right\|}{\left\|A\right\|+\hat{\lambda}}\textbf{e}^{(\lambda-\hat{\lambda})t_{k_i}} &\\
	+\big(\textbf{e}^{\left\|A\right\|d}-\textbf{e}^{-\lambda d}\big)\frac{H_2\left\|cBF\right\|}{\left\|A\right\|+\lambda}\!\!\! &\leq \gamma\textbf{e}^{-\lambda d}
			\end{array}  \label{eq:dineq2}  \right.
\end{align}
where $0<\lambda<\hat{\lambda}$. The expression in \eqref{eq:dineq2} can be used to find the range of values of $d$ that ensure that \eqref{eq:thredel} holds. In other words, we aim to find the maximum admissible delay $\epsilon$. First note that for $d=0$ the expression \eqref{eq:dineq2} reduces to $\beta\leq\gamma$, which is true by design. Since both sides of \eqref{eq:dineq2} are continuous functions of $d$, then, there exist an $\epsilon>0$ such that for $d\in[0,\epsilon)$ the inequality \eqref{eq:dineq2} holds.

An explicit, more conservative, solution $\epsilon>0$ that guarantees \eqref{eq:dineq2} can be found as follows. 
Let $\bar{\lambda}=\hat{\lambda}>\lambda$, 
then, the following two inequalities hold for any $d\geq 0$: 
\begin{align}
\left.
	\begin{array}{l l}
   &\beta\textbf{e}^{\left\|A\right\|d}\! + \!\big(\textbf{e}^{\left\|A\right\|d}\!-\!\textbf{e}^{-\hat{\lambda}d}\big)\frac{H_1\left\|cBF\right\|}{\left\|A\right\|+\hat{\lambda}}\textbf{e}^{(\lambda-\hat{\lambda})t_{k_i}} \\
	&+\big(\textbf{e}^{\left\|A\right\|d}-\textbf{e}^{-\lambda d}\big)\frac{H_2\left\|cBF\right\|}{\left\|A\right\|+\lambda} \\
	\leq\!\!\! & \beta\textbf{e}^{\left\|A\right\|d}\! + \!\big(\textbf{e}^{\left\|A\right\|d}\!-\!\textbf{e}^{-\hat{\lambda}d}\big)\frac{H_1\left\|cBF\right\|}{\left\|A\right\|+\hat{\lambda}}\textbf{e}^{(\lambda-\hat{\lambda})t_{k_i}} \\
	&+\big(\textbf{e}^{\left\|A\right\|d}-\textbf{e}^{-\bar{\lambda} d}\big)\frac{H_2\left\|cBF\right\|}{\left\|A\right\|+\lambda} 
			\end{array}  \label{eq:epseq}  \right.
\end{align}
and
\begin{align}
\left.
	\begin{array}{l l}
\gamma\textbf{e}^{-\bar{\lambda} \tau} \leq  \gamma\textbf{e}^{-\lambda \tau}
			\end{array}  \label{eq:epseq2}  \right.
\end{align}
Then, the solution $\epsilon>0$ of 
\begin{align}
\left.
	\begin{array}{r r}
 \beta\textbf{e}^{\left\|A\right\|d}\! + \!\big(\textbf{e}^{\left\|A\right\|d}\!-\!\textbf{e}^{-\hat{\lambda}d}\big)\frac{H_1\left\|cBF\right\|}{\left\|A\right\|+\hat{\lambda}}\textbf{e}^{(\lambda-\hat{\lambda})t_{k_i}}& \\
	+\big(\textbf{e}^{\left\|A\right\|d}-\textbf{e}^{-\bar{\lambda} d}\big)\frac{H_2\left\|cBF\right\|}{\left\|A\right\|+\lambda}
		&\!\!\!\!= \gamma\textbf{e}^{-\bar{\lambda} \tau}
 \end{array}  \label{eq:epseq3}  \right.
\end{align}
guarantees that inequality \eqref{eq:dineq2} holds. Such solution is given by
\begin{align}
\epsilon=\frac{\ln\Big(\frac{\gamma+H_3}{\beta+ H_3}\Big)}{\left\|A\right\|+\hat{\lambda}}  \label{eq:epsilon}
\end{align}
where
\begin{align}
  H_3=\left\|cBF\right\|\Big(\frac{H_1\textbf{e}^{(\lambda-\hat{\lambda})t_{k_i}}}{\left\|A\right\|+\hat{\lambda}} +\frac{H_2}{\left\|A\right\|+\lambda} \Big)  \label{eq:H3}
\end{align}
By the selection $\hat{\lambda}>\lambda$, we have that $\textbf{e}^{(\lambda-\hat{\lambda})t_{k_i}}\leq 1$ for any $t_{k_i}\geq 0$, and the term $H_3$ remains bounded for any triggering time instant $t_{k_i}>0$, ensuring that $\epsilon>0$. 

Finally, we can show that the minimum inter-event time intervals are strictly positive by following a similar treatment to the one shown in the proof of Theorem \ref{th:consensus}. For this task we need to analyze the response of the error $e_i(t)$ (the error used to trigger events) during the time intervals $t\in[t_{k_i},t_{k_i+1})$. The response of $e_i(t)$ can be bounded as follows 
\begin{align}
\left.
	\begin{array}{l l}
   \left\|e_i(t)\right\|\!\!\!&\leq \int_{t_{k_i}}^t \textbf{e}^{\left\|A\right\|(t-s)}\left\|cBF\right\|\big(H_1\textbf{e}^{-\hat{\lambda}s}+H_2\textbf{e}^{-\lambda s}\big)ds  \\
	&\leq \big(\textbf{e}^{\left\|A\right\|\tau}-\textbf{e}^{-\hat{\lambda}\tau}\big)\frac{H_1\left\|cBF\right\|}{\left\|A\right\|+\hat{\lambda}}\textbf{e}^{-\hat{\lambda}t_{k_i}} \\
	&~~+\big(\textbf{e}^{\left\|A\right\|\tau}-\textbf{e}^{-\lambda\tau}\big)\frac{H_2\left\|cBF\right\|}{\left\|A\right\|+\lambda}\textbf{e}^{-\lambda t_{k_i}}
			\end{array}  \label{eq:eidelb}  \right.
\end{align}
for $t\in[t_{k_i},t_{k_i+1})$, where $\tau=t-t_{k_i}$. Following similar steps to \eqref{eq:taueq}-\eqref{eq:EqTau} we can prove that the inter-event time intervals are lower-bounded by $\tau$, that is, $0<\tau<t_{k_i+1}-t_{k_i}$, where $\tau$ is now given by \eqref{eq:tausoldel}.  $\square$

\section{Example} \label{sec:Examples}
Consider six agents described by \eqref{eq:agents} with
 \[A=
\begin{bmatrix}
	0.192&-0.439\\
	0.431&0.108
\end{bmatrix}, B=
\begin{bmatrix}
	-1.45\\
	0.93
\end{bmatrix}.\]
The solution of the LMI in \eqref{eq:LMI} is given by the following positive definite matrix
\[
P=
\begin{bmatrix}
	0.6174&0.1385\\
	0.1385&0.2754
\end{bmatrix}
\]

The agents are interconnected using a directed communication graph with adjacency matrix given by
 \[\mathcal{A}=
\begin{bmatrix}
	0&0&1&0&0&0\\
	1&0&0&0&0&0\\
	0&1&0&0&0&1\\
	0&0&1&0&0&0\\
	0&1&0&0&0&1\\
	0&1&0&0&0&0
\end{bmatrix}, 
\]
For the system response parameters $\hat{\lambda}=0.24$, and $\hat{\beta}=2$ we select the following parameters: $\lambda=0.03$, $\beta=3$, and $\gamma=12$. Both, $\epsilon$ in \eqref{eq:epsilon} and $\tau$ in \eqref{eq:tausoldel}, depend on the current value of the triggering instant $t_{k_i}$ because the term $H_3$ is a function of this time instants. 

\begin{figure}
	\begin{center}
		\includegraphics[width=8.4cm,height=7cm,trim=1cm .2cm .8cm .5cm]{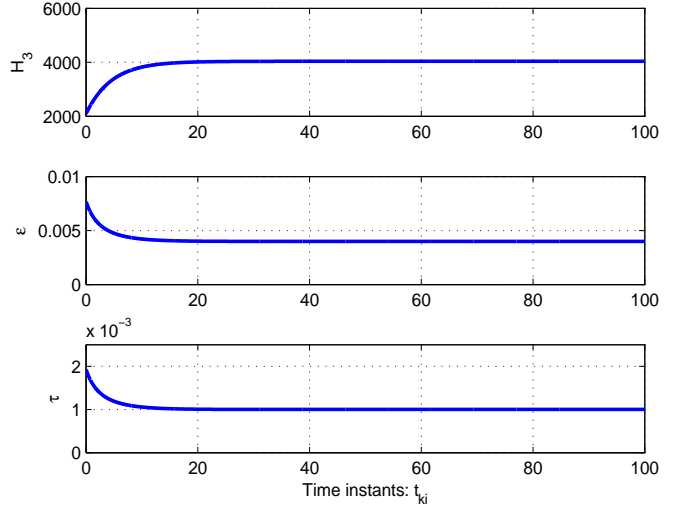}
	\caption{The terms $H_3$, $\epsilon$, and $\tau$ as a function of the triggering time instants $t_{k_i}\geq 0$}
	\label{fig:Het}
	\end{center}
\end{figure}

Fig. \ref{fig:Het} shows the values of these variables for different values of $t_{k_i}\geq 0$. It can be seen that the $H_3\rightarrow\left\|cBF\right\|\frac{H_2}{\left\|A\right\|+\lambda}$ as $t_{k_i}\rightarrow \infty$. Then, we have that $\epsilon$ and $\tau$ converge to constant values around $0.004$ seconds and $0.001$ seconds, respectively.

\begin{figure}
	\begin{center}
		\includegraphics[width=8.4cm,height=7cm,trim=1cm .2cm .8cm .2cm]{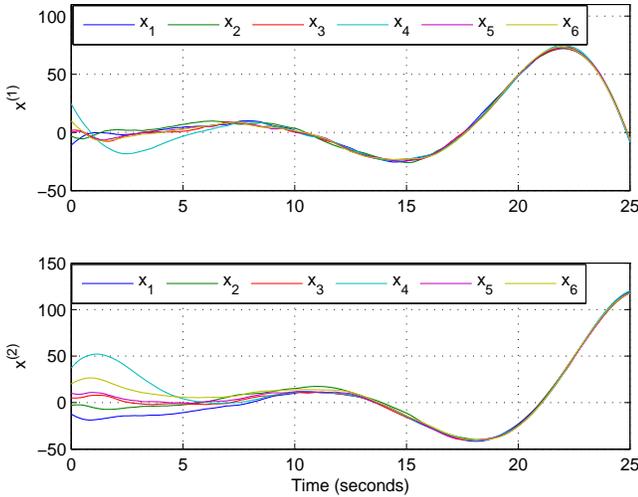}
	\caption{The states of six agents converging to the same trajectories in each dimension}
	\label{fig:states}
	\end{center}
\end{figure}

\begin{figure}
	\begin{center}
		\includegraphics[width=8.4cm,height=7cm,trim=1cm .2cm .8cm .2cm]{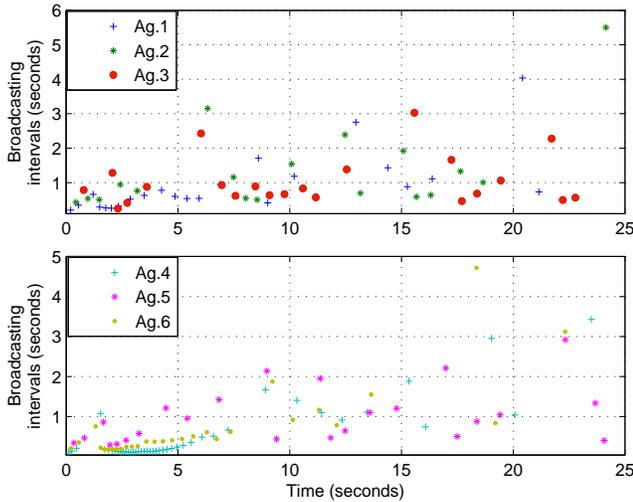}
	\caption{Broadcasting periods for every agent}
	\label{fig:comm}
	\end{center}
\end{figure}

Fig. \ref{fig:states} shows the response of the six agents, for communication delay $d=0.004$ seconds, where each element of the states of the agents converge to the same trajectory. Fig. \ref{fig:comm} shows the transmission periods for every agent where it can be seen that no agent transmits information faster than the lower-bound $\tau=0.001$ seconds.

\section{Conclusions} \label{sec:Conclusions}
Consensus of multi-agent systems with general linear dynamics was discussed in this paper. In this work, agents were not able to communicate continuously and the implementation of decentralized broadcasting strategies was addressed. An event-triggered control technique was presented and it was shown that agents achieve consensus asymptotically. This result applies to the general case where agents' interconnection is represented  by a directed graph. Asymptotic consensus of multi-agent systems under directed graphs and subject to communication delays was also shown in this paper.

\end{document}